\documentclass[preprinst]{elsarticle}\usepackage{amsmath,amssymb,amsfonts,amscd,mathrsfs}
\pdfoutput=1
\usepackage{graphicx,fancybox,color}
\usepackage{arydshln,MnSymbol}
\usepackage{pdfsync,upgreek,dsfont}
\usepackage{enumerate}
\usepackage[utf8]{inputenc}
\renewcommand{\L}{\check{L}}
\newcommand{\A}{{\cal A}}
\newcommand{\B}{{\cal B}}
\newcommand{\C}{{\cal C}}
\newcommand{\D}{{\cal D}}
\newcommand{\X}{{\cal X}}
\renewcommand{\arraystretch}{1}
\newcommand{\ars}{\renewcommand{\arraystretch}{1.2}}
\newcommand{\arc}{\arraycolsep0.5ex}

\newcommand{\skipthis}[1]{}

\renewcommand{\t}[1]{ \tilde{#1} }

\newcommand{\h}[1]{ \widehat{#1} }

\renewcommand{\c}[1]{{\cal #1}}

\newcommand{\epro}{ \hfill\mbox{\rule{2mm}{2mm}} }

\newcommand{\R}{ {\mathbb{R}} }
\newcommand{\N}{ {\mathbb{N}} }

\newcommand{\diag}{ \operatornamewithlimits{diag} }

\renewcommand{\ker}{ \operatorname{ker} }
\newcommand{\im}{ \operatorname{im} }

\renewcommand{\l}[1]{\label{#1}}
\renewcommand{\r}[1]{(\ref{#1})}

\newtheorem{theo}{Theorem}

\newtheorem{hypo}[theo]{Hypothesis}
\newtheorem{lemm}[theo]{Lemma}
\newtheorem{coro}[theo]{Corollary}
\newtheorem{defi}[theo]{Definition}
\newtheorem{rema}[theo]{Remark}

\newcommand{\theorem}[1]{\begin{theo}#1\end{theo}}

\newcommand{\lemma}[1]{\begin{lemm}#1\end{lemm}}

\newcommand{\mul}[1]{\begin{multline}#1\end{multline}}
\newcommand{\equ}[1]{\begin{equation}#1\end{equation}}
\newcommand{\eqn}[1]{$$#1$$}
\newcommand{\eql}[2]{\begin{equation}\label{#1}#2\end{equation}}
\newcommand{\gat}[1]{\begin{gather}#1\end{gather}}

\newcommand{\gan}[1]{\begin{gather*}#1\end{gather*}}

\newcommand{\mas}[2]{\left[\begin{array}{#1}#2\end{array}\right]}
\newcommand{\mat}[2]{\left(\begin{array}{#1}#2\end{array}\right)}

\newcommand{\ite}[1]{\begin{itemize}#1\end{itemize}}

\newcommand{\ga}{\gamma}

\newcommand{\te}[1]{\text{\ \ #1\ \ }}

\newcommand{\cl}{\prec}
\newcommand{\cg}{\succ}

\newcommand{\opt}{{\text{opt}}}

\catcode`\=13
\def{\relax}
\newcommand{\hl}{\\\hline}

\newcommand{\He}{\mbox{He}}
\newcommand{\hdl}{\\ \hdashline}

\renewcommand{\h}{\hat}
\newcommand{\Ga}{\Gamma}

\newcommand{\remark}[1]{\begin{rema}#1\end{rema}}

\usepackage{tikz,pgfplots,datatool}
\usetikzlibrary{positioning,calc,arrows,shapes,shadows}
\tikzset{
auto,
sys/.style 2 args={
rectangle,
draw,
drop shadow,
fill=white,
minimum height=#2,
minimum width=#1,
inner sep=\dn},
sum/.style={circle,draw,draw=black,inner sep=0mm,minimum size=2mm},
jun/.style={circle,draw,draw=black,inner sep=0mm,minimum size=0mm},
>={latex},%triangle 45,
}

\newcommand{\tion}[3]{\coordinate (#1) at ($(#2.north west)!#3!(#2.north east)$)}
\newcommand{\tiow}[3]{\coordinate (#1) at ($(#2.south west)!#3!(#2.north west)$)}
\newcommand{\tioe}[3]{\coordinate (#1) at ($(#2.south east)!#3!(#2.north east)$)}

\def\dn{1ex}
\def\dl{3*\dn}

\tikzstyle{sy0}=[sys={0*\dn}{0*\dn}]
\tikzstyle{sy1}=[sys={12*\dn}{8*\dn}]
\tikzstyle{sy2}=[sys={8*\dn}{6*\dn}]
\tikzstyle{sy3}=[sys={5*\dn}{5*\dn}]
\journal{Systems and Control Letters}
\begin{document}
\allowdisplaybreaks

\begin{frontmatter}
\title{Structured $H_\infty$-Optimal Control for Nested
Interconnections: A State-Space Solution}
\author[carsten]{Carsten W. Scherer\fnref{a1}}
\address[carsten]{Department of Mathematics, University of Stuttgart, \\ Pfaffenwaldring 5a,
70569 Stuttgart, Germany}
\ead{carsten.scherer@mathematik.uni-stuttgart.de}
\fntext[a1]{
The author would like to thank the German Research Foundation (DFG) for financial support of the project within the Cluster of Excellence in Simulation Technology (EXC 310/2) at the University of Stuttgart.
}

\begin{abstract}
If imposing general structural constraints on controllers,
it is unknown how to design $H_\infty$-controllers by convex optimization. Under a so-called quadratic invariance structure of the generalized plant, the Youla parametrization allows to translate the structured synthesis problem into an infinite-dimensional convex program. Nested interconnections that are characterized by a standard plant with a block-triangular structure fall into this class. Recently it has been shown how to design optimal $H_2$-controllers for such nested structures in the state-space by solving algebraic Riccati equations. In the present paper we provide a state-space solution of the corresponding output-feedback $H_\infty$ synthesis problem without any counterpart in the literature. We argue that a solution based on Riccati equations is - even for state-feedback problems - not feasible and we illustrate our results by means of a simple numerical example.
\end{abstract}
\end{frontmatter}

\section*{Notation}
\noindent All matrices in this paper are real and $I$ denotes the identity matrix; in $I_n$ the subscript $n$ specifies its dimension; empty matrices are denoted by $[]$. The range space, the kernel of a matrix $M$ are denoted as $\im(M)$, $\ker(M)$ and $M$ is called a basis matrix of a subspace $\c S\subset\R^n$ if
$M$ is of full column rank (f.c.r.) with $\im(M)=\c S$.  We use the abbreviations $\He(M):=M^T+M$ and $\text{col}(M_1,\ldots,M_n)=(M_1^T\ \ \cdots\ \ M_n^T)^T$. If the transfer matrix $G(s)$ is realized as $C(sI-A)^{-1}B+D$ we write $G=\arc\begin{scriptsize}\mas{c|c}{A&B\hl C&D}\end{scriptsize}$.
Objects that can be inferred by symmetry or are irrelevant are indicated by $\star$.

\section{Introduction}

For $p\in\N$ let us consider the linear time-invariant generalized plant \cite{ZhoDoy96}
\eql{gp}{
\mat{c}{e\hl y_1\\\vdots\\y_p}=
\mat{c|ccc}{
P_{0}   &P_{01} &\cdots &P_{0p} \hl
P_{10}  &P_{11} &\cdots & 0     \\
\vdots  &\vdots &\ddots &\vdots \\
P_{p0}  &P_{p1} &\cdots &P_{pp} }
\mat{c}{d\hl u_1\\\vdots\\u_p}
}
with a control channel that is described by a lower block-triangular transfer matrix of dimension $k\times m$, with a partition into $p$ block-rows and block-columns according to the dimensions
$k=k_1+\cdots+k_p$ and $m=m_1+\cdots+m_p$ respectively.

The goal is to design an internally stabilizing controller
\eql{co}{
\mat{c}{u_1\\\vdots\\u_p}=
\mat{cccc}{
K_{11} &\cdots & 0     \\
\vdots &\ddots &\vdots \\
K_{p1} &\cdots &K_{pp} }
\mat{c}{y_1\\\vdots\\y_p}
}
for \r{gp} which shares the block-triangular structure for the row/column partition
$m=m_1+\cdots+m_p$/$k=k_1+\cdots+k_p$ with the control-channel transfer matrix of \r{gp}
and which renders a bound on the $H_\infty$-norm of the performance channel $d\to e$ satisfied. Let us stress
that the open-loop transfer matrices $d\to e$, $d\to y$ and $u\to e$ do not need to obey any structural constraints.

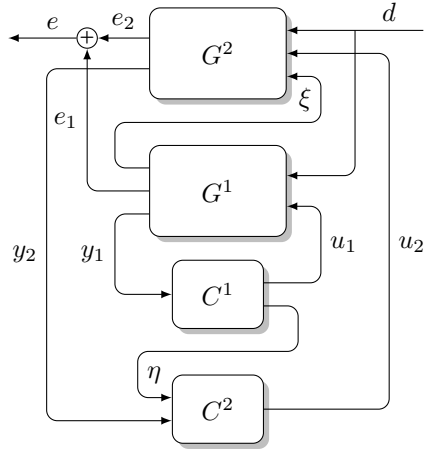
\begin{figure}\center
\begin{tikzpicture}[every path/.style={rounded corners}]
\def\xs{[xshift=-2*\dl]}
\node[sy1] (g) at (0,0)  {$G^1$};
\tiow{go1}{g}{2/4};\tiow{go2}{g}{1/4};\tiow{go3}{g}{3/4};
\tioe{gi1}{g}{2/3};\tioe{gi2}{g}{1/3};\tion{gi3}{g}{2/3};

\node[sy1, above = 4*\dn of g] (h) {$G^2$};
\tiow{ho1}{h}{2/3};\tiow{ho2}{h}{1/3};
\tioe{hi1}{h}{3/4};\tioe{hi2}{h}{2/4};\tioe{hi3}{h}{1/4};

\node[sum, left = 1.5*\dl of ho1] (s) {$+$};

\node[sy2, below = 2*\dn of g] (k) {$C^1$};
\tioe{ko1}{k}{1/3};\tioe{ko2}{k}{2/3};
\node[sy2, below = 4*\dn of k] (l) {$C^2$};
\tiow{li1}{l}{1/3};\tiow{li2}{l}{2/3};

\draw[->] (ko1) -| ([xshift=\dl,yshift=-4*\dn] ko1) -| node[pos=.75]{$\eta$}  ([xshift=-\dl] li2) -- (li2);
\draw[->] (go3) -| ([xshift=-\dl,yshift=4*\dn] go3) -| node[pos=.75]{$\xi$}  ([xshift=\dl] hi3) -- (hi3);

\draw[->] (go1) -|  node[pos=.75]{$e_1$} (s);
\draw[->] ([xshift=2*\dl] hi1)  |- (gi1) ;
\draw[->] (s) --  node[swap,pos=.35]{$e$} ([xshift=-2*\dl] s.west);

\draw[->] (go2)--([xshift=-\dl] go2)|-node[swap,pos=.25]{$y_1$} (k);
\draw[<-] (gi2)--([xshift=\dl] gi2)|-node[pos=.25]{$u_1$} (ko2);

\draw[->] (ho1) --  node[swap]{$e_2$} (s);
\draw[->] ([xshift=4*\dl] hi1)  -- node[swap,pos=.25]{$d$} (hi1) ;
\draw[->] (ho2) --  ([xshift=-3*\dl] ho2) |- node[pos=.2626,swap]{$y_2$}  (li1);
\draw[<-] (hi2) --  ([xshift=3*\dl] hi2)   |- node[pos=.2685]{$u_2$} (l);
\end{tikzpicture}
\caption{A Nested Interconnection}\l{Fig1}
\end{figure}

Such a configuration results e.g. from the nested interconnection in Figure \ref{Fig1} as found and motivated in
\cite{QiSal04}; note that the latter reference provides
other relevant structures of practical interest.

The above formulated synthesis problem has been shown to be
tractable by convex optimization techniques through a structured version of the classical Youla parametrization \cite{GooSer99,Vou00,QiSal04}. In \cite{RotLal06} the authors have shown that quadratic invariance is the essential structural property that allows convexification along this path. The resulting infinite dimensional optimization problem is handled with a Galerkin-type approach, by reducing the search for the structured Youla parameter to a sequence of subspaces of increasing dimension. As its main disadvantage, this approach neither allows
to impose a priori bounds on the degree of (close to) optimal controllers, nor on the to-be-solved optimization problem. For approaching optimality, one might need to rely on high-dimensional subspaces; the incurred numerical instabilities could render it difficult to apply these techniques to large-scale systems.

On the other hand, it is well-known how to tackle the unstructured synthesis problem ($p=1$) in the state-space without the Youla parametrization, either by solving Riccati equations \cite{DoyGlo89} or linear matrix inequalities (LMIs) \cite{GahApk94,IwaSke94}. For the so-called two player problem
($p=2$) and the $H_2$-norm as a cost, such a solution of the structured synthesis problem has been recently proposed in \cite{LesLal13}. However, this approach heavily relies on the inner-product properties of the $H_2$-norm and does, hence, not admit immediate extensions to the $H_\infty$-norm. A similar limitation can be recognized for \cite{ShaPar13}, which handles the $H_2$-problem for more general structures, but is limited to state-feedback synthesis.

Actually, only under strong hypotheses (such as in \cite{Sch02laa}), exact state-space solutions for structured $H_\infty$-synthesis have been available so far. The goal of this paper is to provide
a direct LMI solutions for the described structured $H_\infty$-problem by output-feedback, with a controller construction that is analogous to the one for the unstructured case.
The approach is based on a structured version of the projection lemma as first proposed in \cite{Sch95scl}. If the McMillan degree of \r{gp} is $n$, we show that (almost) optimal controllers have at most degree $np$ (in line with \cite{LesLal13,ShaPar13}), and that they can be constructed by solving a system of LMIs with an a priori fixed dimension.

As a key technical contribution, we show how the LMI framework
allows the design of controllers whose McMillan degree is larger than that of the underlying system, and how to exploit the extra controller dynamics to enforce the required structure. This could pave the way for solutions of other open synthesis problems (e.g.\ in multi-objective control) which remains to be explored. In particular, it could be beneficial in tackling the $H_2$- and $H_\infty$-synthesis problems for general poset structures
\cite{ShaPar13} or even for systems with quadratic invariance \cite{RotLal06}.

The paper is structured as follows. In Section \ref{S2} we review the classical LMI solution of the $H_\infty$-problem and present our new extension for nested systems. We address the construction of controllers and some issues concerning computational complexity. In Section \ref{S3} we provide explicit formulas for the two-block (or two-player) problem and discuss why a solution based on algebraic Riccati equations is out of reach. The paper is concluded with a numerical example in Section \ref{S4}, while the appendices A and B comprise the two technical proofs of this paper.

\section{LMI existence conditions for structured $H_\infty$-controllers}\l{S2}
Let \r{gp} admit the state-space realization
\eql{sys}{
\mat{c}{e\\ y}=
\mas{c|ccc}{
A     &B_{0}  &B \hl
C_{0} &D_0    &E \\
C     &F      &0 }
\mat{c}{d\\ u}
}
in which $A$, $B$ and $C$ share their lower block-triangular structure with the transfer matrix $u\to y$. More specifically, $A\in\R^{n\times n}$ is assumed to be partitioned according to $n=n_1+\cdots+n_p$, which fixes the partition structure of $B\in\R^{n\times m}$ and $C\in\R^{k\times n}$; then all blocks in these matrices above the block-diagonal are supposed to vanish. It is easily seen that such a realization exists \cite{Vou00}. The controller is described similarly as
\eql{con}{
\ars
u=
\mas{c|ccc}{
A^K &B^K\hl
C^K &D^K}y
\te{with lower block-triangular}A^K,\ B^K,\ C^K,\ D^K.
}
The choice of the partition $n^K=n^K_1+\cdots+n^K_p$ of $A^K$ is actually part of the design problem; the ones of $B^K$, $C^K$ and $D^K$ are then determined through $n^K\times k$, $m\times n^K$ and $m\times k$ respectively. Such a controller is said to be structured, in contrast to unstructured controllers that are defined with matrices $A^K$, $B^K$, $C^K$, $D^K$ without any specific sparsity pattern.

For compact notations, we describe triangularly structured matrices $M$ as
\eql{sm1}{\arc
M=\mat{cccc}{
M_{11} &0      &\cdots & 0     \\
M_{21} &M_{21} &\cdots & 0     \\
\vdots &\vdots &\ddots &\vdots \\
M_{p1} &M_{p1} &\cdots &M_{pp} }=
\sum_{j=1}^p L_jM_jR_j^T,\ \
M_j=
\mat{c}{
M_{jj}  \\
M_{2j}  \\
\vdots  \\
M_{pj}  }
\te{for}j=1,\ldots,p}
or through the constraints
\eql{sm2}{\arc
\L_j^TM R_j=0\te{for}j=1,\ldots,p.
}
by making use of the following parts of the block-identity matrix:
$$
R_j=\mat{ccccccc}{0\\\vdots\\0\\I\\\vdots\\0},\ \
L_j=\mat{ccccccc}{0&\cdots&0\\\vdots&&\vdots\\0&\cdots&0\\I&\cdots&0\\\vdots&\ddots&\vdots\\0&\cdots&I}
\te{and}
\L_j=\mat{ccccccc}{I&\cdots&0\\\vdots&\ddots&\vdots\\0&\cdots&I\\0&\cdots&0\\\vdots&&\vdots\\0&\cdots&0}
,\ \ j=1,\ldots,p+1,
$$
Here we highlight the $(j-1)$-st and $j$-th block row for clarity; note that $L_j$ and $\L_j$ have $p-j+1$ and $j-1$ block columns respectively; moreover, we note that $L_1=I$, $\L_{p+1}=I$, $(\L_j\ L_j)=I$
and emphasize $R_{p+1}:=[]$, $L_{p+1}:=[]$, $\L_1:=[]$.
It is stressed that the dimensions of the identity blocks in $R_j$, $L_j$, $\L_j$ can differ (over columns) and are not indicated in the notation; they are determined through the context of the use of these matrices; e.g.\ in \r{sm1} and \r{sm2} those of $\L_j$, $L_j$ and $R_j$ are determined through the row
and column partition of $M$, respectively.

Specifically, the controller matrices in \r{con} are parameterized as
\eql{com}{
\mat{cc}{
D^K  & C^K \\
B^K  & A^K }=
\sum_{j=1}^p
\mat{cc}{L_j&0\\0&L_j}
S_j^K
\mat{cc}{R_j^T&0\\0&R_j^T}
}
with matrices $S_j^K$, $j=1,\ldots,p$, that are unstructured.
It is well-known that the controlled system, the interconnection  of \r{sys} and \r{con}, is described with
\eql{clm}{
\mat{cc}{\A&\B\\\C&\D}:=
\mat{cc|ccc|c}{
A       &0  &B_{0}    \\
0       &0  &0        \hl
C_{0}   &0  &D_{0}    }+
\mat{cccc}{
B     &0     \\
0     &I     \hl
E     &0     \\
}
\mat{cccccccc}{
D^K  & C^K \\
B^K  & A^K }
\mat{cc|cccc|c}{
C     &0   &F   \\
0     &I   &0   }.
}
For some $\ga>0$, the $H_\infty$-design problem consists of finding a controller \r{con} which renders
$\c A$ Hurwitz and such that $\|\C(sI-\A)^{-1}\B+\D\|_\infty<\ga$ is satisfied.
With the classical bounded real lemma \cite{ZhoDoy96}, these two closed-loop properties are equivalently translated into the existence of some
$\X=\X^T$ which satisfies
\eql{lmi1}{
\c X\cg 0,\ \
\He
\mat{ccccc}{X\A&\X\B&0\\0&-\frac{\ga}{2}I&0\\\C&\D&-\frac{\ga}{2}I}=\mat{ccccc}{\A^T\X+\X\A&\X\B&\C^T\\\B^T\X&-\ga I&\D^T\\\C&\D&-\ga I}\cl 0.
}

Let us recall the following LMI conditions for the existence of such controllers that are unstructured
\cite{GahApk94,IwaSke94}. This requires to choose basis matrices $\Phi$ and $\Psi$
of $\ker\mat{ccc}{C&F&0}$ and $\ker\mat{ccc}{B^T&0&E^T}$ respectively.

\theorem{\label{Th1}
There exists an unstructured controller such that the closed-loop system satisfies \r{lmi1} for some symmetric $\c X$ iff there exist symmetric solutions $X$, $Y$ of the following system of LMIs:
\gat{\arc
\label{sy1}
%\mat{ccc}{\Phi_1&0\\ \Phi_2&0\\0&I}^T
\Phi^T
\mat{cccccc}{
A^TX+XA    &XB_0    &C_0^T    \\
B_0^TX     &-\ga I  &D_0^T    \\
C_0        &D_{0}   &-\ga I   }
\Phi\cl0,\ \
%\\\arc
%\label{sy2}
%\mat{ccc}{\Psi_1&0\\ 0&I\\\Psi_2&0}^T
%\ars
\Psi^T
\mat{ccccc|c}{
AY+YA^T   &B_{0}  &YC_0^T                 \\
B_0^T     &-\ga I &D_0^T                 \\
C_{0}Y    &D_{0}  &-\ga I   }
\Psi\cl0,
\\
\label{sy3}
\mat{cc}{Y&I\\I&X}\cg 0.
}
If these LMIs are feasible, one can construct an unstructured controller with McMillan degree at most $n$ which solves the $H_\infty$-problem.
}

Clearly, these conditions have to be incorporated in a problem solution for structured controllers. Let us first describe how the corresponding system of LMIs is composed. Determine basis matrices $\Gamma_j$ with
\eql{bas}{\ars
\ker
\mat{ccccccc|cc}{
\L_j^TC  &\L_j^TF  &0         \\L_j^TB^T & 0       &L_j^TE^T  }=\im(\Ga_j)\te{for}j=1,\ldots,{p+1}. }
Note that we can choose $\Ga_1=\Psi$ and $\Ga_{p+1}=\Phi$. Moreover define
\eql{xy}{\arc
X_j:=\mat{cc}{\h X_{j}&\h Z_{j}^T\\0&I}\te{and}Y_j:=
\mat{cc}{I&0\\-\h Z_{j}&\h Y_{j}}
\te{for}j=1,\ldots,p+1
}
with the decision variables
\eql{dev}{
\h X_j=\h X_j^T,\ \ \h Y_j=\h Y_j^T\te{and unstructured}\h Z_j
}
of dimension $(n_1+\cdots+n_{j-1})$, $(n_j+\cdots+n_{p})$ and
$(n_j+\cdots+n_{p})\times(n_1+\cdots+n_{j-1})$, respectively.
The extreme cases $j=1$, $j={p+1}$ are interpreted as $X_1=I$, $Y_1=Y_1^T$ ($\h X_1$, $\h Z_1$ are empty) as well as $X_{p+1}=X_{p+1}^T$, $Y_{p+1}=I$ ($\h Z_{p+1}$, $\h Y_{p+1}$ are empty).
Note that $X_j^TY_j=\diag(\h X_{j},\h Y_j)$
is symmetric and depends linearly on $\h X_j$, $\h Y_j$. By inspection,
also $X_j^TY_{j+1}$ for $j=1,\ldots,p$ are linear in the decision variables.
Since $L_\nu^TX_j=L_\nu^T$ for $\nu=j,\ldots,p+1$ and $R_\nu^TY_j=R_\nu^T$ for $\nu=1,\ldots,{j-1}$, we infer
$$
X_j^TAY_j=\sum_{\nu=1}^p X_j^T(L_\nu A_\nu R_\nu^T)Y_j=
\sum_{\nu=1}^{j-1} X_j^TL_\nu A_\nu R_\nu^T+
\sum_{\nu=j}^p L_\nu A_\nu R_\nu^TY_j.
$$
Therefore, also $X_j^TAY_j$ depends for $j=1,\ldots,p+1$ linearly on the decision variables. We are now ready to formulate the solution of the $H_\infty$-synthesis problem for structured controllers, the main result of this paper.

\theorem{\label{Th2}
There exists a structured controller such that the closed-loop system satisfies \r{lmi1} with some symmetric $\c X$ iff the following LMIs are feasible:
\gat{
\label{st1}
\Ga_j^T
\mat{cccccc}{
(X_j^TAY_j)^T+X_j^TAY_j    &X_j^TB_0 &Y_j^TC_0^T    \\
B_0^TX_j                   &-\ga I   &D_0^T    \\
C_0Y_j                     &D_{0}    &-\ga I   }
\Ga_j\cl0,\ \ j=1,\ldots,p+1,
\\
\label{st2}\ars
\mat{cc}{Y_j^TX_j&X_j^TY_{j+1}\\Y_{j+1}^TX_{j}&Y_{j+1}^TX_{j+1}}\cg 0\te{for}j=1,\ldots,p.
}
If this LMI system is feasible, one can construct a structured controller with McMillan degree at most $np$ which solves the $H_\infty$-synthesis problem.}

Note that \r{st1}-\r{st2} gracefully specialize to \r{sy1}-\r{sy3} in Theorem \ref{Th1} for $p=1$. The above discussion clarifies that \r{st1} is a system of $p+1$ decoupled LMIs of dimension $n+\dim(d)+\dim(e)$, each defined in the $n(n+1)/2$ scalar decision variables in $(\h X_j,\h Y_j,\h Z_j)$. Furthermore,
\r{st2} involves $p$ LMIs of dimension $2n$, each of which only couples
$(\h X_j,\h Y_j,\h Z_j)$ and $(\h X_{j+1},\h Y_{j+1},\h Z_{j+1})$, respectively. The overall number of scalar decision variables is given by $(p+1)n(n+1)/2$.

If \r{st1}-\r{st2} are feasible, the controller construction proceeds as follows:
\ite{
\item Set $\arc X:=\mat{ccc}{X_2&\cdots&X_{p+1}}$ and $\arc Y:=\mat{ccc}{Y_1&\cdots&Y_{p}}$ and recall $X_1=I$, $Y_{p+1}=I     $.
\item Define $U_1:=0\in\R^{np\times n}$, $\arc U:=\mat{ccc}{U_2&\cdots&U_p}:=I_{np}$, $V_{p+1}:=0$ and
$\arc V:=\mat{ccc}{V_1&\cdots&V_p}$ with
$\arc
V_j:=L_j\mat{ccc}{Y_{j+1}&\cdots&Y_{p+1}}^TX_j-L_j\mat{ccc}{X_{j+1}&\cdots&X_{p+1}}^TY_j
$ for $j=1,\ldots,p.$
\item
Let $K_\nu$ be a basis matrix of
$
\ker
\mat{cc|cccc|c}{
\L_\nu^TCY    &\L_\nu^TC   &\L_\nu^TF &0   \\
\L_\nu^TV     &0           &0         &0   }$, $\nu=1,\ldots,p$.
For $\mu=p,p-1,\ldots,1$, the inequalities
\mul{\label{coco}
\He\left[K_\nu^T\mat{cc|ccc|c}{
AY        &A           &B_{0}             &0                 \\
X^TAY     &X^TA        &X^TB_0            &0                 \hl
0         &0           &-\frac{\ga}{2} I  &0                 \\
C_{0}Y    &C_0         &D_{0}             &-\frac{\ga}{2} I  }K_\nu\right.+\\+
\sum_{j=\nu}^p\left.\arc
K_\nu^T\mat{cccc}{
BL_j    &0       \\
X^TBL_j &U^TL_j  \\
0       &0       \\
EL_j    &0       }S^K_j\!\mat{cc|cccc|c}{
R_j^TCY    &R_j^TC   &R_j^TF &0   \\
R_j^TV     &0        &0      &0   }K_\nu\right]\cl0
}
can be recursively solved algebraically for $S_p^K,S_{p-1}^K\ldots,S_1^K$, by just viewing \r{coco} as an inequality in $S_\nu^K$ and applying the standard projection lemma.
\item With arbitrary solutions of \r{coco}, the controller \r{com} admits the desired block-triangular structure (with $A^K$ partitioned according to $pn=n+\cdots+n$) and solves the $H_\infty$-synthesis problem.
}

If also viewing $\ga$ as a decision variable, one can directly determine the infimal value $\ga_\opt$ of $\ga>0$ for which \r{st1}-\r{st2} are feasible; this optimal value $\ga_\opt$ of the structured
$H_\infty$-synthesis problem can hence be computed through solving an LMI problem of fixed dimension.
Moreover, for any $\ga>\ga_\opt$, one can constructively determine a controller of degree at most $np$ (not depending on $\ga$) that achieves the bound $\ga$.

This is in stark contrast to solutions based on the Youla-parametrization. The possibility to compute $\ga_\opt$ with a semi-definite program of fixed dimension and to guarantee an a priori bound on the degree of suboptimal controllers in the $H_\infty$-setting is considered to be the key novel contribution of this paper.

\remark{\label{rem}
As proved in \ref{A2}, \r{st2} is equivalent to the LMIs
\eql{st2a}{
\ars
\mat{cc|cc}{
\h Y_{j}^{11}    &\h Y_{j}^{12}                    & \h Z_{j}^{1}              &I               \\
\h Y_{j}^{21}    &\h Y_{j}^{22}-\h Y_{j+1}         & \h Z_{j}^{2}-\h Z_{j+1}^1 &-\h Z_{j+1}^2   \hl
(\h Z_{j}^{1})^T &(\h Z_{j}^{2})^T-(\h Z_{j+1}^1)^T& \h X_{j+1}^{11}-\h X_j    &\h X_{j+1}^{12} \\
I                &-(\h Z_{j+1}^2)^T                & \h X_{j+1}^{21}           &\h X_{j+1}^{22} }\cg 0
\te{for}j=1,\ldots,p
}
of the smaller dimension $2n_j+n$, if introducing the more refined partitions
$
\h X_{j+1}=\mat{cc}{\h X_{j+1}^{11}&\h X_{j+1}^{12}\\\h X_{j+1}^{21}&\h X_{j+1}^{22}},\ \
\h Z_{j+1}=\mat{cc}{\h Z_{j+1}^{1}&\h Z_{j+1}^{2}}\te{with}\h X_{j+1}^{22}\in\R^{n_{j}\times n_{j}}
$
and
$
\h Y_{j}=\mat{cc}{\h Y_{j}^{11}&\h Y_{j}^{12}\\\h Y_{j}^{21}&\h Y_{j}^{22}},\ \
\h Z_{j}=\mat{cc}{\h Z_{j}^{1}\\\h Z_{j}^{2}}\te{with}
\h Y_{j}^{11}\in \R^{n_j\times n_j}
$
respectively.}

\section{A specialization and the link to Riccati inequalities/equations}\label{S3}

In view of \cite{LesLal13} and for the purpose of clarity, let us render the two-block case $p=2$ more explicit. Then the matrices in \r{sys} admit the structure
$$
\mat{c|c:cc}{
A     &B_{0}  &B \hl
C_{0} &D_0    &E \hdl
C     &F      &0 }
=
\mat{cc|c:cc}{
A_{11}&0     &B_{1} &B_{11}&0     \\
A_{21}&A_{22}&B_{2} &B_{21}&B_{22}\hl
C_{1}&C_{2}  &D_{0}  &E_{1} &E_{2} \hdl
C_{11}&0     &F_{1}  &0     &0     \\
C_{21}&C_{22}&F_{2}  &0     &0     }.
$$
After dropping indices, the inequalities \r{st1} for $j=1,3$ are equivalent to those in \r{sy1}
in the unstructured symmetric matrices $X$, $Y$ (partitioned according to $A$).
For $j=2$, \r{st1} involves the variables $\h X=\h X^T$, $\h Y=\h Y^T$, $\h Z$ of dimension $\dim(A_{11})$, $\dim(A_{22})$, $\dim(A_{22})\times \dim(A_{11})$ and the basis matrix
$\Ga$ of $\ker\mat{cccc}{C_{11}&0&F_1&0\\0&B_{22}^T&0&E_2^T}$; it reads
explicitly as
\eqn{
\arc
\Ga^T
\mat{ccccccc}{
A_{11}^T\h X+\h XA_{11}       &\!\!(\h Z A_{11}-A_{22}\h Z+A_{21})^T\!\! &\h X B_{1}        &(C_{1}-C_{2}\h Z)^T \\
\h Z A_{11}-A_{22}\h Z+A_{21} &A_{22}\h Y+\h YA_{22}^T          &\h Z B_{1}+B_{2}   &\h Y C_{2}^T        \\
B_{1}^T\h X                   &(\h Z B_{1}+B_{2})^T             &-\ga I             &D_0^T               \\
C_{1}-C_{2}\h Z               &C_{2}\h Y                        &D_{0}              &-\ga I   }
\Ga\cl0.}
In view of Remark \ref{rem}, the coupling conditions \r{st2} can be expressed as
\skipthis{
\eqn{
\mat{cc|cc|ccc}{
Y_{11}  &Y_{12}  &I         & 0          & I        & 0        \\
Y_{21}  &Y_{22}  &-\h Z     & \h Y       & 0        & I        \hl
I       &-\h Z^T &\h X      & 0          & \h X     & 0        \\
0       &\h Y    &0         & \h Y       & \h Z     & I    \hl I       &0       &\h X      & \h Z^T     & X_{11}   & X_{12}   \\
0       &I       &0         & I          & X_{22}   & X_{21}   }\cg 0.
}
Note that this is easily seen to be equivalent to
\eqn{
\mat{cc|cc|ccc}{
Y_{11}  &Y_{12}  &I         & 0          & 0             & 0        \\
Y_{21}  &Y_{22}  &-\h Z     & \h Y       & -\h Z         & I        \hl
I       &-\h Z^T &\h X      & 0          & 0             & 0        \\
0       &\h Y    &0         & \h Y       & \h Z          & I        \hl
I       &\h Z^T  &0         & \h Z^T     & X_{11}-\h X   & X_{12}   \\
0       &I       &0         & I          & X_{22}        & X_{21}   }\cg 0.
}
and hence
\eqn{
\mat{cc|cc|ccc}{
Y_{11}  &Y_{12}       &I         & 0          & 0             & 0        \\
Y_{21}  &Y_{22}-\h Y  &-\h Z     & 0          & 0             & 0        \hl
I       &-\h Z^T      &\h X      & 0          & 0             & 0        \\
0       &0            &0         & \h Y       & \h Z          & I        \hl
0       &0            &0         & \h Z^T     & X_{11}-\h X   & X_{12}   \\
0       &0            &0         & I          & X_{22}        & X_{21}   }\cg 0.
}
}
\eql{st3b}{
\mat{cccc}{
Y_{11}  &Y_{12}       &I          \\
Y_{21}  &Y_{22}-\h Y  &-\h Z      \\
I       &-\h Z^T      &\h X       }\cg 0,\ \
\mat{cccc|ccc}{
 \h Y       & \h Z          & I        \\
 \h Z^T     & X_{11}-\h X   & X_{12}   \\
 I          & X_{22}        & X_{21}   }\cg 0.
}

Let us now explore whether one can decouple these conditions if working with Riccati inequalities or Riccati equations as in \cite{LesLal13}. For simplicity assume
\eql{reg}{\arc
\mat{cc}{B_0\\F}\mat{cc}{B_0\\F}^T=\mat{cc}{\h B_0&0\\0&I}
\te{and}
\mat{cc}{C_0&E}^T\mat{cc}{C_0&E}=\mat{cc}{\h C_0&0\\0&I}.}
Then $\Phi$, $\Psi$, $\Ga$ can be determined explicitly. Elementary computations \cite{GahApk94} show that \r{sy1} and \r{st1} for $j=2$ are equivalent to the following algebraic Riccati inequalities:
\gat{\label{ari1}
A^TX+XA-\ga C^TC+\h C_0/\ga+X\h B_0X/\ga\cl 0,\\
\label{ari2}
AY+YA^T-\ga BB^T+\h B_0/\ga+Y\h C_0Y/\ga\cl 0,\\
\nonumber
\He(X_2^TAY_2)-\ga C^TR_1R_1^TC-\ga BR_2R_2^TB^T+X_2^T\h B_0X_2/\ga+Y_2^T\h C_0Y_2/\ga\cl 0.
}
Because of $R_1^TC=R_1^TCY_2$ and $R_2^TB^T=R_2^TB^TX_2$, the latter one can be transformed by congruence into
\eql{ari3}{
AZ+ZA^T-\ga ZC^TR_1R_1^TCZ-\ga BR_2R_2^TB^T+\h B_0/\ga+Z\h C_0Z/\ga\cl 0}
for the symmetric matrix $Z:=Y_2X_2^{-1}.$
Similarly, \r{st3b} is equivalent to
\eql{st3c}{Y\cg Z\cg 0\te{and}X\cg Z^{-1}.}
The inequalities \r{ari1}-\r{ari2} with the coupling condition $Y\cg X^{-1}\cg 0$ appear in standard $H_\infty$-control \cite{GahApk94} and impose convex constraints on $X$, $Y$. If e.g. $(A,\h B_0)$ is controllable and $(A,\h C_0)$ is observable, one can work with the largest (anti-stabilizing) solutions of the corresponding Riccati equations in order to algebraically  verify the existence conditions for unstructured controllers.

Let us now argue why such a reformulation in terms of Riccati equations seems not possible for structured synthesis. Indeed, both the quadratic and the constant term of the Riccati inequality in \r{ari3} are indefinite. Therefore \r{ari3} imposes a non-convex constraint on $Z$ (and also on $Z^{-1}$). Moreover, feasibility of the inequality does, in general, not imply the existence of a solution of the corresponding Riccati equation. Neither can we easily work with largest or smallest solutions, since the solution set of \r{ari3} does not admit the nice structural properties as known for  \r{ari1}-\r{ari2} \cite{Sch91diss}. It is hence very unclear how to find, among all solutions of the Riccati equation corresponding to \r{ari3}, one which fulfils the coupling condition \r{st3c} as well.
This is yet another motivation for staying with the direct LMI approach as proposed in this paper.

For state-feedback synthesis ($C=I$ and $F=0$) and under the assumptions \r{reg} (without $FF^T=I$),
the existence conditions are easily seen\ to be
\gan{AY+YA^T-\ga BB^T+\h B_0/\ga+Y\h C_0Y/\ga\cl 0,\\
A_{22}\h Y+\h YA_{22}^T-\ga B_{22}B_{22}^T +
(\h Z B_{1}+B_{2})(\h Z B_{1}+B_{2})^T/\ga+
\h Y C_{2}^TC_{2}\h Y/\ga \cl0,\\
\mat{cccc}{
Y_{11}  &Y_{12}         \\
Y_{21}  &Y_{22}-\h Y    }\cg 0\te{and}\h Y\cg 0.
}
Although convex, no reduction to Riccati equations seems possible either.

\section{A numerical example}\label{S4}
Let us perform a simple numerical experiment for the configuration in Fig.~\ref{Fig1} with the stable systems
$$
G^1=\mas{rr|rrr}{-4&2&2&0&0\\1&-0.6&0&0&1\\\hline1&0&0&0&0\\1&1&0&0&1\\0&1&1&0&0}\te{and}
G_{\pm}^2=\mas{rr|rrrr}{-2&1&0&20&0&1\\3&-1.6&0&0&1&0\\\hline1&0&0&0&\pm 0.1\rho&1\\0&1&0&\rho&0&0},
$$
the second one being affected by the parameter $\rho\in[-2,0]$; all signals but $d$ have one component and $d$ has two. In the sequel we calculate the optimal $H_\infty$-levels $\ga_{\rm full}$ for an unstructured controller (with degree $6$ and the results shown in blue), $\ga_{\rm str}$ for a structured controller according to the novel algorithm (with degree $12$ and levels given in red), and $\ga_l$ for a controller designed with the Youla-parametrization as in \cite{QiSal04}. In the latter case,
$l\in\N_0$ indicates the expansion length in the description of the Youla-parameter through
$$\arc
Q(s)=\mat{cc}{q_1&0\\q_3&q_2}\mat{cc}{b_l(s)&0\\0&b_l(s)}\te{with}b_l(s)=
\text{col}\left(1,\frac{1}{s+1},\ldots,\frac{1}{(s+1)^l}\right)
$$
and with free $q_1,q_2,q_3\in\R^{1\times(l+1)}$. Then $Q$ will have the generic degree $2l$ and is designed
with the technique in \cite{Sch00a}, which results in a controller of
degree $6+2l$. For $l\to\infty$ it is known that $\ga_l$ converges to $\ga_{\rm str}$.

\begin{figure}\center
\def\na{exa1}
\begin{tikzpicture}[scale=.8,mark options={scale=.5}]
\begin{axis}
[width=.9\textwidth,height=.5\textwidth,ylabel=Optimal bound $\ga$ for $G^2_-$,xlabel = Parameter $\rho$,grid=both, minor tick num=0]
\pgfplotsinvokeforeach {1,...,2}{
\addplot +[] table[x index=0,y index=#1,header=false] {\na a.tex};}
\pgfplotsset{cycle list shift=-1}
\pgfplotsinvokeforeach {2,...,4}{
\addplot +[green!50!black,mark=triangle] table[x index=0,y index=#1,header=false] {\na b.tex};}
\end{axis}
\end{tikzpicture}
\vspace{-2ex}
\caption{Unstructured (blue), Structured (red), Youla for $l=4,6,8$ (green)}\label{fig2}
\end{figure}
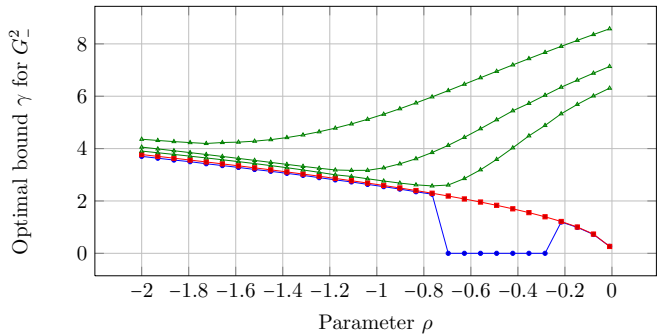

For $G^2_-$ we obtain the results in Fig.~\ref{fig2}. Despite the guaranteed convergence
and $\ga_{\rm full}=\ga_{\rm str}$, the approximation through $\ga_8$ is seen to be of low quality for $\rho\in[-0.2,0]$, while for $\rho\in[-2,0.8]$
the approximation is acceptable, at the expense of a controller order larger than that for the new synthesis technique. We would like to stress that e.g. for $\rho=-2$ and despite the fact that $\ga_{\rm full}=\ga_{\rm str}$, the unstructured controller obtained with Matlab's Robust Control Toolbox
does not ``automatically'' admit the required block-triangular structure.

\begin{figure}[h]\center
\def\na{exa2}
\begin{tikzpicture}[scale=.8,mark options={scale=.5}]
\begin{axis}
[width=.9\textwidth,height=.5\textwidth,ylabel=Optimal bound $\ga$ for $G^2_+$,xlabel = Parameter $\rho$,grid=both, minor tick num=0]
\pgfplotsinvokeforeach {1,...,2}{
\addplot +[] table[x index=0,y index=#1,header=false] {\na a.tex};}
\pgfplotsinvokeforeach {2,...,4}{
\addplot [green!50!black,mark=triangle] table[x index=0,y index=#1,header=false] {\na b.tex};}
\end{axis}
\end{tikzpicture}
\vspace{-2ex}
\caption{Unstructured (blue), Structured (red), Youla for $l=4,6,8$ (green)}\label{fig3}
\end{figure}
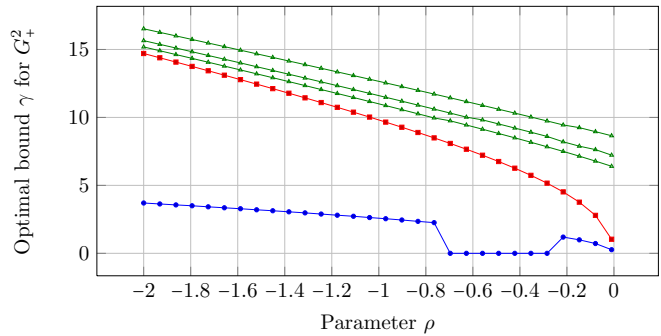

For $G^2_+$ (Fig.~\ref{fig3}) we obtain an example with $\ga_{\rm str}>\ga_{\rm full}$ and with a widening gap for decreasing values of $\rho$. The approximation error $\ga_l-\ga_{\rm str}>0$ is getting again smaller for smaller values of $\rho$. The computation times for the red curves are similar to those for the full controller (and for an example of this size), while the computation of $\ga_8$ takes in our implementation about ten times longer.

\section{Conclusions}\label{S5}

For the first time we have given a direct and exact LMI solution for the optimal $H_\infty$-design of controllers with a block-triangular structure for generalized plants whose control channel match this structure. As the key novel ingredient, both the dimension of the to-be-solved convex optimization problem and the degree of the controller can be fixed a priori in terms of dimensional parameters of the problem data. The proposed technique for designing controllers whose order is larger than that of the generalized plant lends itsself for various generalizations that are currently under investigation.

\section*{References}
%\bibliographystyle{elsart-num-sort}
%\bibliography{x:/Carsten/Literatur/ref/ref}

\appendix
\section{Proof of Theorem \ref{Th2}}\label{A1}

\subsection{Preparation}

We first sketch the proof for the unstructured case. Suppose that a full
controller has been found which renders \r{lmi1} satisfied. By adding uncontrollable (or unobservable)
stable modes in the controller, we can assume w.l.o.g. that $n^K\geq n$. In a partition of
$\c X$ and $\c X^{-1}$
according to $n+n^K$, let us denote by $\mat{c}{X\\U}$ and $\mat{c}{Y\\V}$  the first block columns of
these matrices respectively. Since tall, we can assume w.l.o.g. (by perturbing $\c X$) that $U$ is of full column rank. We directly infer the relations
\eql{fac}{
\c X\c Y=\c Z\te{for}\c Y:=\mat{cc}{Y&I\\V&0},\ \ \c Z:=\mat{cc}{I&X\\0&U}.
}
Since $\c Z$ has full column rank, the same holds for $\c Y$. With a congruence transformation
involving $\c Y$, \r{lmi1} implies
\eql{lmi2}{
\c Y^T\c X\c Y\cg 0,\ \
\He
\mat{ccccc}{
\c Z^T\A\c Y&\c Z^T \B&0\\0&-\frac{\ga}{2}I&0\\\C \c Y&\D&-\frac{\ga}{2}I}\cl 0.}
The first inequality can be expressed as
\eql{sy3p}{
\mat{cc}{Y^T&Y^TX+V^TU\\I&X}\cg 0
}
and is equivalent to \r{sy3} by the symmetry of the left-hand side.
By inspection, the second inequality in \r{lmi2} reads with \r{clm} as
\eql{lmi3}{\arc
\He\left[Q(X,Y)+
\mat{cc|cc}{
B^T   &B^TX &0 &E^T  \\
0     &U    &0 &0    }^T
\mat{cccccccc}{
D^K  & C^K \\
B^K  & A^K }
\mat{cc|cccc|c}{
CY    &C   &F &0   \\
V     &0   &0 &0   }\right]\cl 0,
}
where we introduced the abbreviation
\eql{defQ}{\arc
Q(X,Y):=\mat{cc|ccc|c}{
AY      &A     &B_{0}           &0                 \\
X^TAY   &X^TA  &X^TB_0          &0                 \hl
0       &0     &-\frac{\ga}{2}I &0                 \\
C_{0}Y  &C_0   &D_{0}           &-\frac{\ga}{2}I   }.
}
The key step is the elimination of the controller parameters in
\r{lmi3} with the projection lemma \cite{GahApk94}. With suitable row-partitions of
$\Phi$ and $\Psi$ we have
$$\arc
\ker\mat{cc|cccc|c}{
CY    &C   &F &0   \\
V     &0   &0 &0   }=
\im(\Phi_e)\te{and}
\ker
\mat{cc|cc}{
B^T   &B^TX &0 &E^T  \\
0     &U    &0 &0    }=\im(\Psi_e)
$$
for
$\Phi_e:=\text{col}(0,\Phi_1,\Phi_2,\Phi_3)$ and $\Psi_e:=\text{col}(\Psi_1,0,\Psi_2,\Psi_3)$.
We exploited that $U$ and $V$ have full column rank such that the kernels have
basis matrices that are independent of  $Y$, $V$, $X$, $U$!
Obviously, $\Phi_e^T\r{lmi3}\Phi_e$ and $\Psi_e^T\r{lmi3}\Psi_e$ simplify to
\eql{sy1p}{
\He\left[\Phi_e^TQ(X,Y)\Phi_e\right]=
\He\left[
\Phi^T
\mat{c|ccc|c}{
X^TA  &X^TB_0          &0                 \hl
0     &-\frac{\ga}{2}I &0                 \\
C_0   &D_{0}           &-\frac{\ga}{2}I   }\Phi\right]\cl 0,
}
\eql{sy2p}{
\He
\left[\Psi_e^TQ(X,Y)\Psi_e\right]=
\He
\left[
\Psi^T
\mat{c|cccc|c}{
AY      &B_{0}           &0                 \hl
0       &-\frac{\ga}{2}I &0                 \\
C_{0}Y  &D_{0}           &-\frac{\ga}{2}I   }\Psi
\right]\cl0,
}
which shows that $X$ and $Y$ satisfy \r{sy1}.

Now suppose that \r{sy1}-\r{sy3} hold. Then set $U=I$ and $V=I-YX$ to
make sure that \r{sy3} implies \r{sy3p} since the left-hand sides are identical. Due to \r{sy1p}-\r{sy2p}, there exist
$A^K$, $B^K$, $C^K$, $D^K$ with \r{lmi3}. This nontrivial part of the projection lemma is
constructive and leads to the controller parameters. We then get back to \r{lmi2}. Since $\c Y$
is square and invertible, we can transform \r{lmi2} into \r{lmi1} which proves that the controller does the
required job.

For structured controller synthesis, it is a natural idea to try applying an extension of
the projection lemma that allows for structured unknowns. We will reveal that the following old
generalization from \cite{Sch95scl} serves this purpose.

\lemma{\label{L1}Let us be given $Q$ and $M_1,\ldots,M_{p},M_{p+1}=0$, $N_0=0,N_1,\ldots,N_{p}$,
with $\ker(M_j)\subset\ker(M_{j+1})$ for $j=1,\ldots,p-1$. Then there exists $S_1,\ldots,S_p$ with
$\He\left[Q+\sum_{j=1}^p M_j^TS_jN_j\right]\cl 0$
iff $\He(Q)$ is negative definite on the subspaces
$\ker(N_0)\cap\cdots\cap\ker(N_{j-1})\cap \ker(M_j)\te{for}j=1,\ldots,{p+1}.$
}
If $K_\nu$ denotes a basis matrix of $\ker(N_0)\cap\cdots\cap\ker(N_{\nu-1})$, the proof shows that one can
compute $S_p,S_{p-1},\ldots,S_1$ by solving
$\He\left[K_\nu^T\left(Q+\sum_{j=\nu}^p M_j^TS_jN_j\right)K_\nu\right]\cl 0$
for $S_\nu$ and for $\nu=p,p-1,\ldots,1$ recursively.

\subsection{Proof of Necessity in Theorem \ref{Th2}}

Let there exist a structured controller \r{con} such that \r{lmi1} holds for some $\c X=\c X^T$.
As in the unstructured case we assume w.l.o.g. that $n^K_j\geq n$ for $j=1,\ldots,p$.
With $L_j$, $\L_j$, $R_j$ chosen according to the partition $n_1^K+\cdots+n_p^K$ of $A^K$, define
$$\ars\arc
\mat{ccccc}{\t X&\t U_j^T\\\t U_j&\t Z_j}:=
\mat{cc}{I_n&0\\0&L_j^T}\c X\mat{cc}{I_n&0\\0&L_j}\cg 0,\ \
\mat{ccccc}{\t Y_j\\\t V_j}:=\mat{ccccc}{\t X&\t U_j^T\\\t U_j&\t Z_j}^{-1}\mat{c}{I_n\\0}
$$
for $j=1,\ldots,p+1$.
Clearly $\t U_j\in\R^{(n_j+\cdots+n_p)\times n}$ just consists of the blocks below $\t X$ in the first
block column $\t U_1$ of $\c X$. Since $\t U_j$ is tall and $R_j^TL_j\t Z_j^{-1}$ (which is the first block row of $\t Z_j^{-1}$) has full row rank, we can
slightly perturb the blocks of $\t U_1$ in $\c X$ without violating \r{lmi1} and such that
$
R_j^TL_j\t Z_j^{-1}\t U_j\te{has f.c.r. for all}j=1,\ldots,p.
$
By the block-inversion formula we have
$\t V_j=-\t Z_j^{-1}\t U_j(\t X-\t U_j^T\t Z_j^{-1}\t U_j)^{-1}$ and hence
\eql{frtv}{
R_j^TL_j\t V_j\te{has f.c.r. for}j=1,\ldots,p.}

Let us momentarily zoom in and choose $\L_j$, $L_j$ according to the partition $n_1+\cdots+n_p$ of $A$; then, trivially, $\arc \t Y_j=\small\mat{cc}{\L_j^T\t Y_j\L_j&\L_j^T\t Y_jL_j\\L_j^T\t Y_j\L_j&L_j^T\t Y_jL_j}\cg 0$; with the blocks
$$
\h X_j:=(\L_j^T\t Y_j\L_j)^{-1},\
\h Y_j=L_j^T\t Y_jL_j-(L_j^T\t Y_j\L_j)\h X_j(\L_j^T\t Y_jL_j),\
\h Z_j^T=-\h X_j(\L_j^T\t Y_jL_j)
$$
we define the invertible matrices $X_j$ and $Y_j$ by \r{xy} in order to get $\t Y_jX_j=Y_j$.

It follows directly from the definitions that
\eql{eq1}{
\mat{cc}{I_n&0\\0&L_j^T}\c X\mat{ccccc}{\t Y_j\\L_j\t V_j}=
\mat{ccccc}{I_n\\0}\te{for}j=1,\ldots,p+1.
}
By right-multiplying \r{eq1} with $X_j$ and if setting $V_j:=L_j\t V_jX_j$, we get
\eql{eq2}{\mat{cc}{I_n&0\\0&L_j^T}\c X\mat{ccccc}{Y_j\\V_j}=\mat{ccccc}{X_j\\0}\te{for}j=1,\ldots,p+1.}
For suitable $U_j$ we hence arrive at the equation
\eql{facs}{
\c X
\mat{ccccc}{
Y_1  &Y_2  &\cdots&Y_p  &Y_{p+1}\\
V_1  &V_2  &\cdots&V_p  &V_{p+1}}=
\mat{cccccc}{
X_1  &X_2  &\cdots&X_p  &X_{p+1} \\
U_1  &U_2  &\cdots&U_p  &U_{p+1} }.
}
Due to \r{eq2}, $\L_j^TL_j=0$ and \r{frtv} together with $\det(X_j)\neq 0$, we infer that
\eql{pro}{\arc
L_j^TU_j=0,\ \ \L_j^TV_j=0,\ \ R_j^TV_j\te{has f.c.r. for}j=1,\ldots,p+1.
}
Since $L_\nu^TU_\nu=0$ and $\ker(L_\nu^T)\subset\ker(L_j^T)$ for $\nu=1,\ldots,j$, we can actually conclude $L_j^TU_\nu=0$ for $\nu\leq j$; similarly $\L_j^TV_\nu=0$ for $\nu\geq j$. Hence \r{pro} actually implies
\eql{zer}{\arc
L_j^T\mat{ccc}{U_1&\cdots&U_j}=0\te{and}\L_j^T\mat{ccc}{V_j&\cdots&V_{p+1}}=0
\te{for}j=1,\ldots,p+1,
}
Note that $V_{p+1}=0$ and $U_1=0$. In fact, the relations \r{zer} mean that
\eql{vu}{\arc V:=\mat{ccc}{V_1&\cdots&V_p}\te{and}U:=\mat{ccc}{U_2&\cdots&U_{p+1}}}
are lower and upper block-triangular matrices, respectively. This is the reason why the rank conditions in \r{pro} even imply that
\eql{frv}{\arc
\L_{j}^T\mat{ccc}{V_1&\cdots&V_{j-1}}\te{and}L_{j}^T\mat{ccc}{V_j&\cdots&V_p}
\te{have f.c.r. for}j=1,\ldots,p+1.
}
Indeed, if $x=\text{col}(x_1,\ldots,x_{j-1})$ is in the kernel of the first matrix, we infer for
$\mu=2,\ldots,j$ with $\ker(\L_\mu^T)\supset\ker(\L_j^T)$ that
$0=\L_{\mu}^T\sum_{\nu=1}^{j-1} V_\nu x_\nu=\L_{\mu}^T\sum_{\nu=1}^{\mu-1} V_\nu x_\nu$
and hence $0=R_{\mu-1}^T\sum_{\nu=1}^{\mu-1} V_\nu x_\nu$.
With $\mu=2$ we get $R_1^TV_1 x_1=0$, i.e. $x_1=0$; for $\mu=3$ we
thus obtain $R_2^TV_2 x_2=0$, i.e. $x_2=0$; similarly one shows $x_3=0,\ldots,x_{j-1}=0$, i.e.  $x=0$. Analogous arguments apply for the second matrix in \r{frv}.

We can finally prove that
\eql{fru}{\arc
L_{j}^T\mat{ccc}{U_{j+1}&\cdots&U_{p+1}}\te{has f.c.r. for}j=1,\ldots,p;}
if we exploit $\L_j\L_j^T+L_jL_j^T=I$ and \r{zer}, the relation \r{facs} implies
$$\arc\ars
\mat{cc}{I_n&0\\0&L_j^T}
\c X
\mat{cc}{I_n&0\\0&L_j}
\mat{ccccc}{
Y_{j}         &\cdots &Y_p     &Y_{p+1}\\
L_j^TV_{j}    &\cdots &L_j^TV_p&0      }=
\mat{cccccc}{
X_{j}       &X_{j+1}       &\cdots&X_{p+1} \\
0           &L_j^TU_{j+1}  &\cdots&L_j^TU_{p+1} }.
$$
Since $\c X$ is invertible and $L_j$, $Y_{p+1}$, $X_j$  have f.c.r., the claim follows
from \r{frv}.

Let us now recall $X_1=I_n$ and $Y_{p+1}=I_n$. Therefore, \r{facs} can be compactly expressed as
\r{fac} with \r{vu}, $Y:=(Y_1\ \ \cdots\ \ Y_p)$ and $X:=(X_{2}\ \ \cdots\ \ X_{p+1})$.
By \r{frv}, $V$ has f.c.r. and hence the same holds for $\c Y$.
As for unstructured controllers, \r{lmi2} hence implies \r{lmi3}. In view of \r{com}, this reads as
\eql{stp}{
\arc
\He\left[
Q(X,Y)\!+\!\!\sum_{j=1}^p
\mat{cccc}{
L_j^TB^T &L_j^TB^T\!X &0   &L_j^TE^T  \\
0        &L_j^TU      &0   &0         }^T\!\!\!\!S^K_j\!
\mat{cc|cccc|c}{
R_j^TCY    &R_j^TC   &R_j^TF &0   \\
R_j^TV     &0        &0      &0   }\right]\!\cl\!0.
}
Now we apply Lemma \ref{L1}. This is possible since the kernels $\ker(L_j^T)$ and hence also
$
\ker
\mat{cccc}{
L_j^TB^T &L_j^TB^TX &0   &L_j^TE^T  \\
0        &L_j^TU      &0   &0       }
$
form a non-decreasing sequence of subspaces. Let us also observe that
$$\ars
\bigcap_{\nu=1}^{j-1}
\ker
\mat{cc|cccc|c}{
R_\nu^TCY    &R_\nu^TC   &R_\nu^TF &0   \\
R_\nu^TV     &0        &0      &0   }=
\ker
\mat{cc|cccc|c}{
\L_j^TCY    &\L_j^TC   &\L_j^TF &0   \\
\L_j^TV     &0        &0      &0   }.
$$
To use Lemma \ref{L1}, we need to determine basis matrices of the intersections of the latter two kernels, which just equals the kernel of
\equ{\label{ker}
\ars\arc
\mat{ccccccc|cc}{
\L_j^TCY_1  &\!\cdots\! &\L_{j}^TCY_{j-1}  &\L_j^TCY_j &\L_{j}^TCY_{j+1} &\!\cdots\!  &\L_j^TCY_{p+1}  &\L_j^TF  &0         \\
\L_j^TV_1   &\!\cdots\! &\L_{j}^TV_{j-1}   &\L_j^TV_j  &\L_{j}^TV_{j+1}  &\!\cdots\!  &\L_l^TV_{p+1}   &0        &0         \hl
L_j^TB^TX_1 &\!\cdots\! &L_{j}^TB^TX_{j-1} &L_j^TB^TX_j&L_{j}^TB^TX_{j+1}&\!\cdots\!  &L_j^TB^TX_{p+1} & 0       &\!\!\!L_j^TE^T  \\
L_j^TU_1    &\!\cdots\! &L_j^TU_{j-1}      &L_j^TU_{j} &L_{j}^TU_{j+1}   &\!\cdots\!  &L_j^TU_{p+1}    & 0       &0         }.
}
At this point the relevance of the triangular structure of $U$, $V$ and the specific construction of $X_j$, $Y_j$, comes to light; we also crucially exploit that $B$ and $C$ are block-triangular through \r{sm2}.
In view of $\L_\nu^TCR_\nu=0$ and $\ker(\L_j^T)\supset\ker(\L_\nu^T)$ for $\nu=j,\ldots,p+1$,
we actually have $\L_j^TCR_\nu=0$ for $\nu\geq j$ and hence $\L_j^TCL_j=0$.
With $(\L_j\ L_j)(\L_j\ L_j)^T=I$ and $\L_j^TY_\nu=\L_j^T$ for $\nu\geq j$ (due to \r{xy}), we conclude
$$\arc\ars
\L_j^TCY_\nu=
\mat{cc}{\L_j^TC\L_j&\L_j^TCL_j}\mat{cc}{\L_j^TY_\nu\\L_j^TY_\nu}=
\mat{cc}{\L_j^TC\L_j&\L_j^TCL_j}\mat{cc}{\L_j^T\\L_j^T}=
\L_j^TC
$$
for $\nu\geq j$. In exactly the same fashion one proves $L_j^TB^TX_\nu=L_j^TB^T$ for $\nu\leq j$.
If combining with \r{zer}, the matrix \r{ker} actually simplifies to
\equ{\label{ker2}
\ars\arc
\mat{ccccccc|cc}{
\L_j^TCY_1      &\!\cdots\! &\L_{j}^TCY_{j-1}  &\L_j^TC &\L_{j}^TC        &\!\cdots\!  &\L_j^TC  &\L_j^TF  &0                   \\
\L_j^TV_1       &\!\cdots\! &\L_{j}^TV_{j-1}   &0       &0                &\!\cdots\!  &0               &0        &0                   \hl
L_j^TB^T        &\!\cdots\! &L_{j}^TB^T        &L_j^TB^T&L_{j}^TB^TX_{j+1}&\!\cdots\!  &L_j^TB^TX_{p+1} & 0       &L_j^TE^T  \\
0               &\!\cdots\! &0                 &0       &L_{j}^TU_{j+1}   &\!\cdots\!  &L_j^TU_{p+1}    & 0       &0                   }.
}
In view of \r{frv} and \r{fru}, the kernel of \r{ker2} is {\em exactly} given by
\eql{ann}{
\im(\Ga_e)\te{with}\Ga_e:=\text{col}(0,\ldots,0,\Ga_{1j},0,\ldots,0,\Ga_{2j},\Ga_{3j})
}
and with a suitable row partition of $\Ga_j=\text{col}(\Ga_{1j},\Ga_{2j},\Ga_{3j})$. Hence Lemma \ref{L1}
applied to \r{stp} leads with \r{defQ} to the inequalities
\eql{sol}{\arc
0\cg \He\left[\Ga_e^TQ(X,Y)\Ga_e\right]=
\He\left[
\Ga_j^T
\mat{cc|ccc|c}{
X_j^TAY_j   &X_j^TB_0        &0                 \hl
0           &-\frac{\ga}{2}I &0                 \\
C_{0}Y_j    &D_{0}           &-\frac{\ga}{2}I   }\Ga_j\right],\ j=1,\ldots,p+1,
}
which is clearly identical to \r{st1}. Moreover, $\c Y^T\c X\c Y\cg 0$ is identical to
\eql{co1}{\ars\arc
\mat{cc}{Y^T&Y^TX+V^TU\\I&X}=
\mat{c|cccc|c}{
Y_1^T   &Y_1^TX_2+V_1^TU_2&\cdots & Y_1^TX_{p+1}+V_1^TU_{p+1}   \\
\vdots  &                                                       \\
Y_p^T   &Y_p^TX_2+V_p^TU_2&\cdots & Y_p^TX_{p+1}+V_p^TU_{p+1}   \hl
I       &X_2 &\cdots&X_{p+1}
}\cg 0.
}
Due to \r{pro}, $\im(V_j)\subset\ker(\L_j^T)=\im(L_j)\subset\im(L_\nu)$ and $L_\nu^TU_\nu=0$ for $\nu=1,\ldots,j$ imply $V_{j}^TU_\nu=0$. Hence the terms $V_j^TU_\nu$ for $j\geq \nu$ on and below the block-diagonal of the left-hand side of \r{co1} vanish; by symmetry, this inequality hence just reads as
\eql{co2}{
(Y_j^TX_k)_{j,k=1,\ldots,p+1}\cg 0.
}
In view of $X_1=I$ and $Y_{p+1}=I$ this clearly implies \r{st2}.\epro

\subsection{Proof of Sufficiency in Theorem \ref{Th2}}

Suppose \r{st1}-\r{st2} are feasible and set $X_1:=I_n$, $Y_{p+1}:=I_n$. Let us first show that \r{st2}
implies \r{co2}. Since $Y_j^TX_j\cg 0$ for $j=1,\ldots,p+1$, both $X_j$ and $Y_j$ are invertible. By subtracting the $X_{j+1}^{-1}X_j$-right-multiple of column $j+1$ from column $j$ and the $X_j^TX_{j+1}^{-T}$ left-multiple of row $j+1$ from row $j$ (which is a congruence transformation), \r{co2} transforms due to its structure into
\eql{co3}{
\diag_{j=1,\ldots,p}\left(Y_j^TX_j-X_j^TY_{j+1}X_{j+1}^{-1}X_j\right)\cg0.
}
By taking the Schur-complement, this is indeed implied by \r{st2}.

Now choose $n^K_j:=n$ and $U_{j+1}=R_j$ for $j=1,\ldots,p$ as well as $U_{1}:=0$.
Then $U:=(U_2\ \cdots\ U_{p+1})=I_{np}$ and the first relation in \r{pro} is trivially valid.
Set $V_{p+1}:=0$ and choose the block-columns of $V=(V_1\ \cdots\ V_p)$ such that all properties in \r{pro} hold and, in addition, \r{co2} is identical to \r{co1}; though crucial, this is actually easy to establish: Just take
\eql{defV}{\arc
V_j:=L_j\check V_j\te{with } \check V_j^T:=X_j^T\mat{ccc}{Y_{j+1}&\cdots&Y_{p+1}}-Y_j^T\mat{ccc}{X_{j+1}&\cdots&X_{p+1}}.
}
Then $\L_jV_j=0$. Also note that $R_j^TV_j=Y_{j+1}^TX_j-X_{j+1}^TY_j$ which is the negative transpose of $(Y_j^TX_j-X_j^TY_{j+1}X_{j+1}^{-1}X_j)X_j^{-1}X_{j+1}$ and thus invertible by \r{co3}. Therefore, \r{pro} is valid. For $j=1,\ldots,p$, this shows $V_j^T(U_1\ \ \cdots\ \ U_j)=0$; on the other hand, due to $L_j^T(U_{j+1}\ \cdots\ U_{p+1})=L_j^TL_j=I$, we conclude from \r{defV} that
$$\arc
Y_j^T\mat{ccc}{X_{j+1}&\cdots&X_{p+1}}+V_j^T \mat{ccc}{U_{j+1}&\cdots&U_{p+1}}=X_j^T\mat{ccc}{Y_{j+1}&\cdots&Y_{p+1}}.
$$
Hence \r{co2} is identical to \r{co1}.
Moreover, as seen in the necessity proof, \r{zer}, \r{frv} and \r{fru} are valid, which in turn implies again that the kernel of \r{ker} equals \r{ann}.
Now note that \r{st1} is the same as \r{sol}, and that the latter are the precise conditions for the existence of $S^K_1,\ldots,S^K_p$ which satisfy \r{stp} (Lemma \ref{L1}).

We finally claim that the resulting controller \r{com} solves the $H_\infty$-problem. For this purpose we just note that, by construction, $U$ and $V$ are square and invertible. Hence the same holds for $\c Y$, $\c Z$ in \r{fac}, which allows us to define $\c X:=\c Z\c Y^{-1}$. Then \r{co1} and \r{stp} are nothing but \r{sy3p} and \r{lmi3}, respectively.
Since $\c Y$ is invertible, this leads back to \r{lmi2} which proves the claim. The partition of $A^K$ is determined by the block-structure $pn=n+\dots+n$ of $U$. \epro

\section{Proof of equivalence in Remark \ref{rem}}\label{A2}

With the refined partitions, \r{st2} reads as
$$\ars
\mat{ccc|ccc}{
\h X_j  &0                &0                  & \h X_{j}                   & 0              &0          \\
0       & \h Y_{j}^{11}   &\h Y_{j}^{12}      & \h Z_{j}^{1}               &I               &0          \\
0       & \h Y_{j}^{21}   &\h Y_{j}^{22}      & \h Z_{j}^{2}-\h Z_{j+1}^1  &-\h Z_{j+1}^2   &\h Y_{j+1} \hl
\h X_j  &(\h Z_{j}^{1})^T &\star              & \h X_{j+1}^{11}            &\h X_{j+1}^{12} & 0         \\
0       &I                &-(\h Z_{j+1}^2)^T  & \h X_{j+1}^{21}            &\h X_{j+1}^{22} & 0         \\
0       &0                &\h Y_{j+1}         & 0                          & 0              &           \h Y_{j+1}}\cg 0,\ \ j=1,\ldots,p.
$$
By congruence, this is equivalent to \r{st2a} and $\h X_j\cg 0$, $\h Y_{j+1}\cg 0$ for $j=1,\ldots,p$. The latter two inequalities are redundant: We have $\h X_1=I\cg 0$ and $\h Y_{p+1}=I\cg 0$; if assuming
$\h X_j\cg 0$, $\h Y_{j+1}\cg0$, we can directly infer from the right-lower and left-upper block of \r{st2a} that $\h X_{j+1}\cg 0$, $\h Y_j\cg 0$; the statement follows by induction.

\end{document}